# AN EXAMINATION OF THE SPILLAGE DISTRIBUTION

BEN O'NEILL[*], *Australian National University*[**]

WRITTEN 29 DECEMBER 2021

**Abstract**

We examine a family of discrete probability distributions that describes the "spillage number" in the extended balls-in-bins model. The spillage number is defined as the number of balls that occupy their bins minus the total number of occupied bins. This probability distribution can be characterised as a normed version of the expansion of the noncentral Stirling numbers of the second kind in terms of the central Stirling numbers of the second kind. Alternatively, it can be derived in a natural way from the extended balls-in-bins model. We derive the generating functions for this distribution and important moments of the distribution. We also derive an algorithm for recursive computation of the mass values for the distribution. Finally, we examine the asymptotic behaviour of the spillage distribution and the performance of an approximation to the distribution.

EXTENDED BALLS-IN-BINS MODEL; SPILLAGE DISTRIBUTION; NON-CENTRAL STIRLING NUMBERS OF THE SECOND KIND.

Occupancy distributions are an interesting class of distributions that arise in statistical problems involving simple-random-sampling-with-replacement. O'Neill (2021) examines three families of distributions that arise in an extended balls-in-bins model that generalises this case. In the present paper we will examine an interesting discrete distribution that arises in a balls-in-bins model, that is based around a well-known expansion for non-central Stirling numbers of the second kind. Let $S(n, k, \phi)$ denote the non-central Stirling numbers of the second kind and let $S(n, k)$ denote the (central) Stirling numbers of the second kind. These two types of numbers are related by the following expansion (see e.g., Charalambides 2005, p. 77, Eqn. 2.16):[1]

$$S(n, k, \phi) = \sum_{r=0}^{n-k} \binom{n}{k+r} \cdot \phi^{n-k-r} \cdot S(k+r, k).$$

For any values $\phi \geq 0$ for the non-centrality parameter, all terms in this sum are non-negative, which means we can use this sum to form a discrete probability distribution which we will call the spillage distribution (see the definition below). This distribution was first introduced in O'Neill (2021) in the context of a broader examination of probability distributions arising from occupancy distributions in the extended balls-in-bins model, but it also arises naturally as a normed version of the above mathematical expansion.

---

[*] E-mail address: ben.oneill@anu.edu.au; ben.oneill@hotmail.com.
[**] Research School of Population Health, Australian National University, Canberra ACT 0200.

[1] The result we state here is a slight variation of Eqn. 2.16 of Charalambides (2005) (p. 77) where we have used the change-of-variable $r = t - k$. It is worth noting that there is also another similar expansion where we use falling factorials instead of powers in the expansion (Eqn. 2.15, p. 77).

**DEFINITION (The spillage distribution):** This distribution is a discrete probability distribution with probability mass function given by:[2]

$$\text{Spillage}(r|n, k, \phi) \equiv \binom{n}{k+r} \cdot \phi^{n-k-r} \cdot \frac{S(k+r, k)}{S(n, k, \phi)} \qquad r = 0, \ldots, n-k,$$

where $n \in \mathbb{N}$ is the **size parameter**, $0 \leq k \leq n$ is the **occupancy parameter** and $0 \leq \phi \leq \infty$ is the **scale parameter**. □

The spillage distribution arises in the extended balls-in-bins model from the occupancy number and the effective number of balls (see O'Neill 2021). Suppose we randomly allocate $n$ balls to $m$ bins and a ball "occupies" its allocated bin with probability $0 \leq \theta \leq 1$ (and with probability $1 - \theta$ it "falls through" its allocated bin). The **effective balls** $0 \leq n_{\text{eff}} \leq n$ is the total number of balls that occupy their bins and the **occupancy number** $0 \leq K_n \leq \min(n, m)$ is the total number of occupied bins (bins that are occupied by at least one ball). The **spillage number** is the value $R_n = n_{\text{eff}} - K_n$, which can be considered as the number of effective balls beyond those that increase the occupancy number. O'Neill (2021) shows that the spillage distribution is the conditional distribution of the spillage number given the occupancy number:

$$\mathbb{P}(R_n = r | K_n = k, n, m, \theta) = \text{Spillage}\left(r \Big| n, k, m \cdot \frac{1-\theta}{\theta}\right).$$

**REMARK:** This conditional distribution for the spillage number depends on the parameters $m$ and $\theta$ only through the scale parameter $\phi = m \cdot (1 - \theta)/\theta$. In the classical case where $\theta = 1$ all the balls occupy their bins; this means we have $\phi = 0$ and $n_{\text{eff}} = n$ with probability one, giving $R_n = n - k$ with probability one. In the case where $m = \infty$ and $0 < \theta < 1$ all balls are allocated to different bins; this means we have $\phi = \infty$ and $n_{\text{eff}} = k$ with probability one, giving a $R_n = 0$ with probability one. In the case where $m = \infty$ and $\theta = 1$ we have an indeterminate form for the scale parameter. This corresponds to allocation of balls to an infinite number of bins, with zero probability of falling through the bins, so it should give $n_{\text{eff}} = K_n = n$ with probability one; by convention we set $\phi \equiv 0$ in this case so that the distribution for the spillage number is a point mass on the outcome $R_n = n - k$, but this gives strange results when $n \neq k$, since we condition on an occupancy number with zero probability.[3] □

---

[2] The case where $\phi = 0$ gives the point-mass distribution $\text{Spillage}(r|n, k, 0) = \mathbb{I}(r = n - k)$. Additionally, the case where $\phi = \infty$ is defined via the limit $\text{Spillage}(r|n, k, \infty) \equiv \lim_{\phi \to \infty} \text{Spillage}(r|n, k, \phi) = \mathbb{I}(r = 0)$.

[3] In discrete problems, the conditional distribution of a quantity —conditional on an impossible result— can be set to any distribution without causing problems. Regardless, the law-of-total-probability recovers the appropriate marginal distributions for the quantities.

# 1. Properties of the spillage distribution

Here we will derive the generating functions and important central moments of the spillage distribution. The distribution has a closed-form expression for its cumulant function, and the derivatives of this function, which means we can obtain closed-form expressions for all the moments. We begin by establishing the generating functions for the distribution.

**THEOREM 1 (Generating functions):** The spillage distribution has the following generating functions. The probability generating function $G$, characteristic function $\phi$, moment generating function $m$, and cumulant function $K$ are given respectively by:

$$G(z) = z^{n-k} \cdot \frac{S(n, k, \phi/z)}{S(n, k, \phi)},$$

$$\phi(t) = e^{it(n-k)} \cdot \frac{S(n, k, \phi e^{-it})}{S(n, k, \phi)},$$

$$m(t) = e^{t(n-k)} \cdot \frac{S(n, k, \phi e^{-t})}{S(n, k, \phi)},$$

$$K(s) = s(n - k) + \log S(n, k, \phi e^{-s}) - \log S(n, k, \phi).$$

As can be seen from Theorem 1, the generating functions for the spillage distribution can be written in a succinct form, but this form involves ratios of values of the non-central Stirling numbers of the second kind (or logarithms of these ratios). Unfortunately, these functions are just as complicated as the mass function for the distribution. Now, to obtain the moments of the spillage distribution, it is useful to be able to find the derivatives of the cumulant function and its Maclaurin representation. To facilitate this analysis, we define the quantities:

$$Q_\ell(s) \equiv \frac{S(n - \ell, k, \phi e^{-s})}{S(n, k, \phi e^{-s})} \qquad H_\ell \equiv \phi^\ell \cdot \frac{S(n - \ell, k, \phi)}{S(n, k, \phi)}.$$

These quantities appear as terms in the Maclaurin representation of the cumulant function and the expansion for the moment generating function (for the latter, see below). In Lemma 1-2 we establish some preliminary results involving these quantities and then in Theorems 2-3 we establish important central moments for the spillage distribution and their asymptotic forms.

**LEMMA 1:** We have:

$$\frac{dQ_\ell}{ds}(s) = \phi e^{-\ell s}[nQ_\ell(s)Q_1(s) - (n - \ell)Q_{\ell+1}(s)].$$

**LEMMA 2:** The moment generating function for the spillage distribution can be written as:

$$m(t) = e^{t(n-k)} \sum_{\ell=0}^{n-k} \binom{n}{\ell} (e^{-t} - 1)^\ell \cdot H_\ell .$$

**THEOREM 2 (Moments):** The moments of the spillage distribution are:

$$\mu_{n,k,\phi} \equiv \mathbb{E}(R_n | K_n = k) = (n - k) - nH_1,$$

$$\sigma^2_{n,k,\phi} \equiv \mathbb{V}(R_n | K_n = k) = nH_1 - n^2 H_1^2 + n(n-1)H_2,$$

$$\gamma_{n,k,\phi} \equiv \mathbb{Skew}(R_n | K_n = k) = \frac{\begin{bmatrix} -nH_1 + 3n^2 H_1^2 - 3n(n-1)H_2 - 2n^3 H_1^3 \\ + 3n^2(n-1)H_1 H_2 - n(n-1)(n-2)H_3 \end{bmatrix}}{[nH_1 - n^2 H_1^2 + n(n-1)H_2]^{3/2}},$$

$$\kappa_{n,k,\phi} \equiv \mathbb{Kurt}(R_n | K_n = k) = 3 + \frac{\begin{bmatrix} nH_1 - 7n^2 H_1^2 + 7n(n-1)H_2 + 12n^3 H_1^3 \\ -18n^2(n-1)H_1 H_2 + 6n(n-1)(n-2)H_3 \\ -6n^4 H_1^4 + 12n^3(n-1)H_1^2 H_2 - 3n^2(n-1)^2 H_2^2 \\ -4n^2(n-1)(n-2)H_1 H_3 \\ +n(n-1)(n-2)(n-3)H_4 \end{bmatrix}}{[nH_1 - n^2 H_1^2 + n(n-1)H_2]^2}.$$

**THEOREM 3 (Asymptotic moments):** As $n \to \infty$ we have:

$$H_\ell \to \psi^\ell \qquad \psi \equiv \frac{\phi}{k + \phi}.$$

This gives the corresponding asymptotic moments:

$$\mu_{n,k,\phi} \sim \mu^*_{n,k,\phi} \equiv (n-k)(1-\psi) \qquad \sigma^2_{n,k,\phi} \sim \sigma^{*2}_{n,k,\phi} \equiv (n-k)\psi(1-\psi),$$

$$\gamma_{n,k,\phi} \sim \gamma^*_{n,k,\phi} \equiv \frac{2(\psi - \tfrac{1}{2})}{\sqrt{(n-k)\psi(1-\psi)}} \qquad \kappa_{n,k,\phi} \sim \kappa^*_{n,k,\phi} \equiv 3 + \frac{1 - 6\psi(1-\psi)}{(n-k)\psi(1-\psi)}.$$

Theorem 2 gives the exact form for the main moments of the spillage distribution. Since the moments involve the noncentral Stirling numbers of the second kind they are roughly as difficult to compute as the values in the mass function, and it is usual to compute them using recursive methods. Theorem 3 gives asymptotic forms for the moments which hold when we take $n \to \infty$ (taking the parameters $k$ and $\phi$ to be fixed). From these asymptotic forms we can see that the distribution is asymptotically unskewed and mesokurtic. It is possible to obtain higher-order moments for the distribution, but the general form is extremely cumbersome. Since the moments involve the noncentral Stirling numbers of the second kind they are roughly as difficult to compute as the original mass function, and it is usual to compute them using recursive methods. The asymptotic moments shown in Theorem 3 are identical to the binomial distribution with $n - k$ trials and probability parameter $1 - \psi$. In fact, the spillage distribution converges to the binomial as $n - k$ becomes large (in the sense that these distributions become asymptotically equivalent as $n \to \infty$). And of course, since the binomial is asymptotically equivalent to the normal in this limit, the spillage distribution is also asymptotically equivalent to the normal distribution as $n \to \infty$.

Aside from its generating functions and moments, it is useful to look at a mixture result that involves the spillage distribution. In Theorem 4 below we show that the binomial distribution can be formed as a mixture of spillage distributions, with the occupancy parameter generated by the occupancy distribution. If we take an occupancy number $K \sim \text{Occ}(n, m, \theta)$ and then a spillage number $R|K \sim \text{Spillage}(n, K, m(1-\theta)/\theta)$ then the sum of these numbers (which is the effective number of balls in the model) has distribution $K + R \sim \text{Bin}(n, \theta)$. This result reflects the underlying logic of the occupancy number and spillage number in the extended balls-in-bins model described in O'Neill (2021).

**THEOREM 4 (Occupancy mixture):** The binomial distribution satisfies the equation:

$$\text{Bin}(x|n, \theta) = \sum_{k=0}^{x} \text{Spillage}\left(x - k \middle| n, k, m \cdot \frac{1-\theta}{\theta}\right) \cdot \text{Occ}(k|n, m, \theta).$$

O'Neill (2021) derives several mixture results involving the occupancy distributions appearing in the extended balls-in-bins model. In particular, the paper shows that the extended occupancy distribution can be written as a binomial mixture of classical occupancy distributions as:

$$\text{Occ}(k|n, m, \theta) = \sum_{n_*=k}^{n} \text{Occ}(k|n_*, m) \cdot \text{Bin}(n_*|n, \theta)$$

$$= \sum_{r_*=0}^{n-k} \text{Occ}(k|k + r_*, m) \cdot \text{Bin}(k + r_*|n, \theta).$$

The mass function for the spillage distribution can therefore be written as:

$$\begin{aligned}\text{Spillage}(r|n, k, \phi) &= \mathbb{P}(n_* - k = r | K_n = k) \\ &= \frac{\mathbb{P}(n_* = k + r, K_n = k)}{\mathbb{P}(K_n = k)} \\ &= \frac{\text{Occ}(k|k + r, m) \cdot \text{Bin}(k + r|n, \theta)}{\text{Occ}(k|n, m, \theta)}.\end{aligned}$$

$$= \frac{\text{Occ}(k|k+r,m) \cdot \text{Bin}(k+r|n,\theta)}{\sum_{r_*=0}^{n-k} \text{Occ}(k|k+r_*,m) \cdot \text{Bin}(k+r_*|n,\theta)}.$$

As can be seen, the properties of the spillage distribution tie into the properties of the occupancy distribution and binomial distribution.

## 2. Computation and approximation of the spillage distribution

Computation of the mass function for the spillage distribution is made difficult by the presence of the Stirling numbers of the second kind. It is generally not a good idea to attempt to compute these numbers directly from their explicit form due to arithmetic underflow problems. For this reason it is best to compute these numbers recursively. For the spillage distribution we need to compute the Stirling numbers $S(k,k), \ldots, S(n,k)$ (or equivalent transformations of these values) but this effectively requires us to compute a matrix of all values up to the latter indices. To avoid arithmetic underflow we undertake all computations in log-space, so we give the recursive equation in log-space. In Theorem 5 below we show a useful recursive equation in log-space that can be used to compute the mass function.

**THEOREM 5 (Recursive equation):** Consider the function:

$$L(r,k) = \log \binom{n}{k+r} + (n-k-r)\log\phi + \log S(k+r,k).$$

For all $r \geq 1$ and $k \geq 1$ this function obeys the recursive equations:

$$L(0,0) = n\log\phi,$$
$$L(r,0) = -\infty,$$
$$L(0,k) = \log\binom{n}{k} + (n-k)\log\phi,$$
$$L(r,k) = \log(n-k-r+1) - \log(k+r) - \log(\phi)$$
$$\quad + \text{logsumexp}(\log(k) + L(r-1,k), L(r,k-1)).$$

We then have $\log \text{Spillage}(r|n,k,\phi) = L(r,k) - \text{logsumexp}(L(0,k), \ldots, L(n-k,k))$.

In Algorithm 1 below we show how to compute the mass function using the recursive equations in Theorem 5. The algorithm treats the cases $\phi = 0$ and $\phi = \infty$ as special cases (where the spillage distribution reduces to a point-mass distribution) and uses the recursive algorithm for the remaining cases. The algorithm fixes the values $n$ and $\phi$ and forms an $(r+1) \times (k+1)$ matrix $\mathbf{L}$ containing the values $L(r,k) = \mathbf{L}_{r,k}$. We compute the first row of this matrix using the base equation and compute the remaining rows recursively. The final column of the matrix gives the log-probabilities for the spillage distribution with the occupancy parameter $k$. The algorithm returns either the log-probabilities or probabilities of the spillage distribution.

---

**ALGORITHM 1: Spillage Distribution**

---

**Input:**     Size parameter **n** (positive integer)
              Occupancy parameter **k** (non-negative integer no greater than **n**)
              Scale parameter **φ** (non-negative real value or infinity)
              Logical value **log** (specifying whether output is log-probability)
**Output:**    Vector of probabilities/log-probabilities
              from the spillage mass function over arguments **r = 0,…,n-k**

---

```
#Deal with special case where φ = 0
if (φ = 0) {
  LSPILLAGE <- [-Inf, … , -Inf, 0] indexed by r = 0,…,n-k }

#Deal with special case where φ = Inf
if (φ = Inf) {
  LSPILLAGE <- [0, -Inf, … , -Inf] indexed by r = 0,…,n-k }

#Deal with remaining cases
if (0 < φ < Inf) {

  L <- Matrix with rows 0,…,r and columns 0,…,k
        (Initial values all set to -Inf)
  L[0,0] <- n*log(φ)
  for each kk = 1,…,k {
    L[0,kk] <- log(choose(n,kk)) - (n-kk)*log(φ)
    for each rr = 1,…,r {
      L[rr,kk] <- log(n-kk-rr+1) - log(kk+rr) - log(φ) +
                 logsumexp(log(kk) + L[rr-1,kk], L[rr, kk-1]) } }
  LSPILLAGE <- L[,k] - logsumexp(L[,k]) }

#Return output
if (log) { LSPILLAGE } else { exp(LSPILLAGE) }
```

Algorithm 1 computes the mass function over the range $r = 0, \ldots, n - k$. It can be modified to return the mass values over an arbitrary vector of argument values for $r$. The computational intensity is proportionate to the size of the occupancy parameter $k$ since the algorithm computes values in the log-matrix for values of the occupancy parameter $1, \ldots, k$. Since only the last column of the matrix is used in the mass function, these other intermediate computations are generally wasted. It is possible to modify the algorithm to use the entire block of values and compute "blocks" of mass values for the distribution for all occupancy parameters up to the value of the space parameter. In this case we compute an $(r + 1) \times (\min(n, m) + 1)$ matrix **L** and we use each column of the matrix to obtain the corresponding mass values for the spillage distributions with occupancy parameters $k = 0, \ldots, \min(n, m)$.

Algorithm 1 is computationally feasible so long as the parameters $n$ and $k$ are not too high. However, in cases where both these parameters are high it may be cumbersome to use recursive computation. In this case it is useful to be able to approximate the spillage distribution by an appropriate asymptotic form. Approximating the spillage distribution is relatively simple — it requires us to compute approximations to the Stirling numbers of the second kind, to give the following general form for the approximating distribution:

$$\widehat{\text{Spillage}}(r|n,k,\phi) \propto \binom{n}{k+r} \cdot \phi^{n-k-r} \cdot \hat{S}(k+r,k).$$

There are a substantial number of known approximations to the Stirling numbers of the second kind (Hsu 1948, Moser and Wyman 1958, Good 1961, Harper 1967, Bender 1973, Bleick and Wang 1974, Temme 1993, Louchard 2013). A useful overview of the available approximations in given in Louchard (2013). Most of these approximations are obtained by representing the Stirling numbers in their integral form and then using various saddle-point approximations to the integral. Every convergent approximation to the Stirling numbers of the second kind gives a corresponding approximation for the spillage distribution, so we will not attempt to give anything like an exhaustive treatment of the available options.

To give just one example of an approximation to the distribution, we will use a relatively simple approximation for the Stirling numbers in Bender (1973, pp. 108-109). In this paper it is shown that the Stirling numbers of the second kind are well approximated by the quantities:

$$\hat{S}(k+r,k) = \frac{(k+r)_r}{\sqrt{2\pi(k+r)}} \cdot \frac{e^{-\alpha k}}{[1-e^{\alpha}\log1\text{pexp}(-\alpha)][\log1\text{pexp}(-\alpha)]^{k+r-2}},$$

where $\alpha \equiv \alpha(r,k)$ is defined implicitly by the equation $1+r/k = (1+e^{\alpha})\log1\text{pexp}(-\alpha)$. Substituting this into the approximating kernel of the spillage distribution, and removing any multiplicative terms that do not depend on $r$, gives the kernel:

$$\widehat{\text{Spillage}}(r|n,k,\phi) \propto \binom{n}{k+r} \cdot \phi^{n-k-r} \cdot \hat{S}(k+r,k)$$

$$\propto \binom{n}{k+r} \cdot \frac{(k+r)_r}{\sqrt{k+r}} \cdot \frac{\phi^{n-k-r} \cdot e^{-\alpha k}}{[1-e^{\alpha}\log1\text{pexp}(-\alpha)][\log1\text{pexp}(-\alpha)]^{k+r-2}}$$

$$\propto \frac{1}{(n-k-r)!} \cdot \frac{1}{\sqrt{k+r}} \cdot \frac{\phi^{n-k-r} \cdot e^{-\alpha}}{[\log1\text{pexp}(-\alpha)-r/k][\log1\text{pexp}(-\alpha)]^{k+r-2}}$$

Computation of this approximation proceeds by first computing each of the kernel values and then normalising this kernel to obtain the approximate mass function of the distribution. The computation of the kernel values require us to find $\alpha \equiv \alpha(r,k)$ over $r = 0, \ldots, n-k$, where each of these computations requires us to compute using the implicit function for this value.[4] In Algorithm 2 we compute this approximating distribution in log-space and then convert back to probability-space in the final step.

---

**ALGORITHM 2: Spillage Distribution (Approx)**
**(Based on approximation in Bender 1973)**

**Input:**  Size parameter **n** (positive integer)
Occupancy parameter **k** (non-negative integer no greater than **n**)
Scale parameter **ϕ** (non-negative real value or infinity)
Logical value **log** (specifying whether output is log-probability)
**Output:** Vector of probabilities/log-probabilities
from the spillage mass function over arguments **r = 0,…,n-k**

```
#Set log-kernel vector
LOGKERNEL <- Vector with elements r = 0,…,n-k
            (Initial values all set to -Inf)

#Compute log-kernel vector
for each r = 0,…,n-k {

  #Find the parameter ALPHA
  OBJECTIVE <- function(a) { (1 + r/k - (1+exp(a))*log1pexp(-a))^2 }
  ALPHA     <- Argument that minimises OBJECTIVE

  #Compute the log-kernel
  LOGKERNEL[r] <- - log(k+r)/2 - logfactorial(n-k-r) - ALPHA*k
                 + (n-k-r)*log(ϕ) - (k+r-2)*log(log1pexp(-ALPHA))
                 - log(log1pexp(-ALPHA) - r/k) }

#Compute the approximate log-mass function
APPROX <- LOGKERNEL - logsumexp(LOGKERNEL)

#Return output
if (log) { APPROX } else { exp(APPROX) }
```

---

Algorithm 1 gives "exact" computation of the mass function for the spillage distribution based on recursive computation of the Stirling numbers of the second kind. (It is not actually exact; there is arithmetic error owing to floating-point representation of real numbers.) Algorithm 2 gives a useful approximation to the mass function of the spillage distribution in the case where the parameters $n$ and $k$ are large. In Figure 1 below we show an example of the mass function of the spillage distribution, with the approximation shown in the same figure. The black points show the "exact" values of the mass function (computed using the recursive algorithm) and the

---

[4] There are several possible ways to perform this intermediate computation, but we recommend using optimisation routines to compute $\alpha = \mathrm{argmax}_a F(a)$ using the objective function:

$$F(a) \equiv \left[1 + \frac{r}{k} - (1 + e^a)\log 1\mathrm{pexp}(-a)\right]^2 \qquad a \in \mathbb{R}.$$

red points and lines show the approximate values. We can see that the approximation is already quite good even with the moderate values $n = 100$, $k = 30$ and $\phi = 40$, though there is some systematic deviation from the true values.[5] The accuracy of the approximation improves as we increase these parameter values and the exact and approximate values converge in the limit.

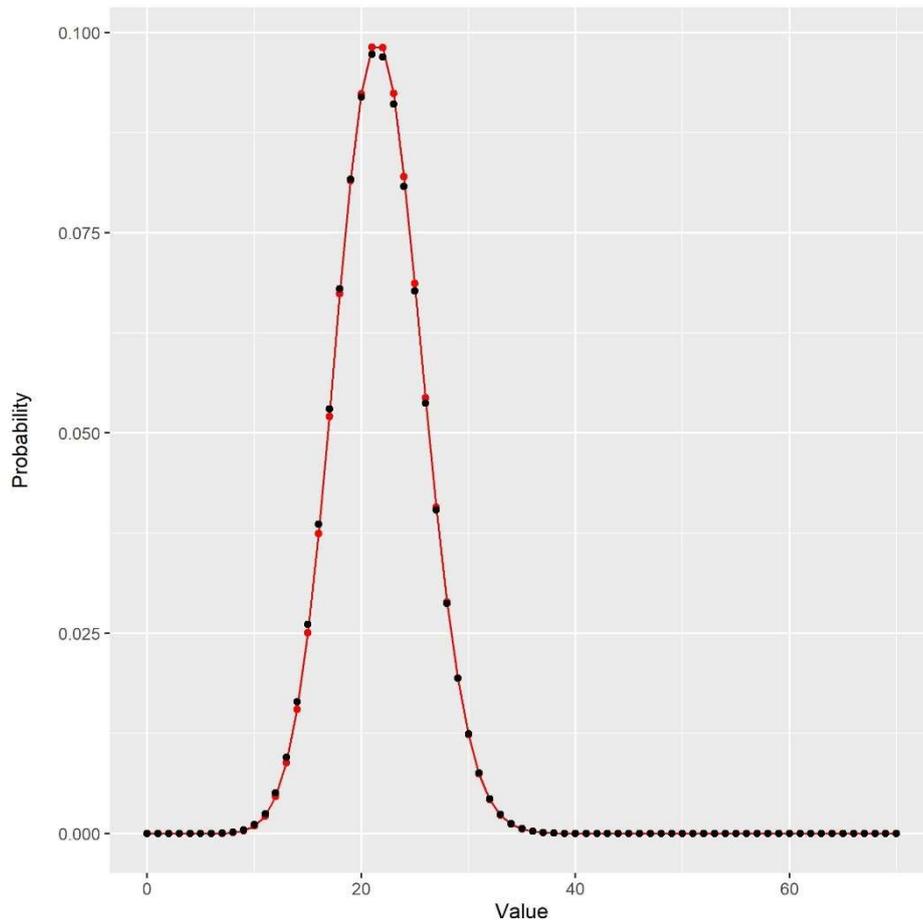

**FIGURE 1:** The spillage distribution and approximation to the distribution

Spillage distribution has parameters $n = 100$, $k = 30$ and $\phi = 40$.
Approximation is the red dots and line; the true mass values are the black dots.

It is simple to create a hybrid algorithm that allows the user to switch between the recursive computation of the mass function or the approximation. It is also relatively simple to extend these algorithms to compute the cumulative distribution function and quantile function for the distribution. These functions for the spillage distribution are implemented in the **occupancy** package in **R** (O'Neill 2021). This package provides standard probability functions for the distribution with an allowance for the user to choose exact or approximate computation.

---

[5] The maximum absolute difference here is $\max_r |\text{Spillage}(r|n,k,\phi) - \widehat{\text{Spillage}}(r|n,k,\phi)| = 0.001294$.

## 3. Accuracy of the approximation

In order to get a better idea of the accuracy of the approximation to the spillage distribution we can compute the log-root-mean-square-error (LRMSE) of the approximate mass function over a range of parameter values. Our analysis is similar to analysis in O'Neill (2020) that computed the LRMSE between the true distribution and its approximation to assess the accuracy of an approximation to the occupancy distribution over a large block of parameter values. In the present case the LRMSE of the approximation is given by:

$$\text{LRMSE}(n, k, \phi) \equiv \frac{1}{2} \cdot \log\left(\frac{1}{n - k + 1} \sum_{r=0}^{n-k} (\text{Spillage}(r|n, k, \phi) - \widehat{\text{Spillage}}(r|n, k, \phi))^2\right).$$

We computed the LRMSE over a range of parameter values using $n = 10, 20, \ldots, 10000$. In Figure 2 we show the accuracy over a range of values up to $n = 1000$. (The grey squares represent perfect accuracy of the approximation, which occurs when the distribution is a point-mass.) As can be seen, the approximation becomes more accurate as $n$ gets larger, but it is also more accurate when the variance of the distribution is larger. For large $n$ the variance of the spillage distribution is maximised when $k \approx \phi$, which is reflected in regions of high accuracy.

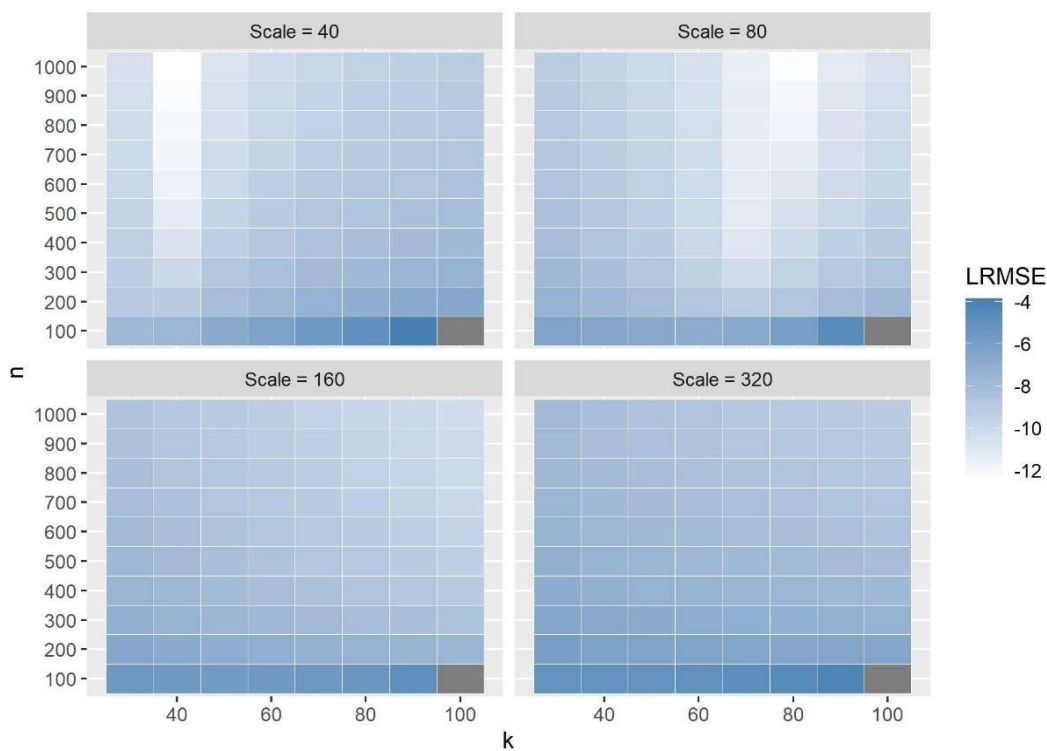

**FIGURE 2:** LRMSE of approximation to spillage distribution

The main thing that affects the accuracy of the approximation to the spillage distribution is the variance of the distribution. Higher variances generally correspond to higher accuracy (lower RMSE) of the approximation. This relationship is exhibited in Figure 3 below, where we show the RMSE for the approximation against the asymptotic variance for the spillage distribution, for each combination of parameter values in our simulation (with each variable on a log-scale).[6] As can be seen from the figure, there is a close relationship between these variables. In fact, with the present simulations, the (asymptotic) log-variance $\sigma^{*2}_{n,k,\phi} = (n-k)\psi(1-\psi)$ and the proportion $\psi$ explains almost all the variance in the accuracy of the approximation.[7]

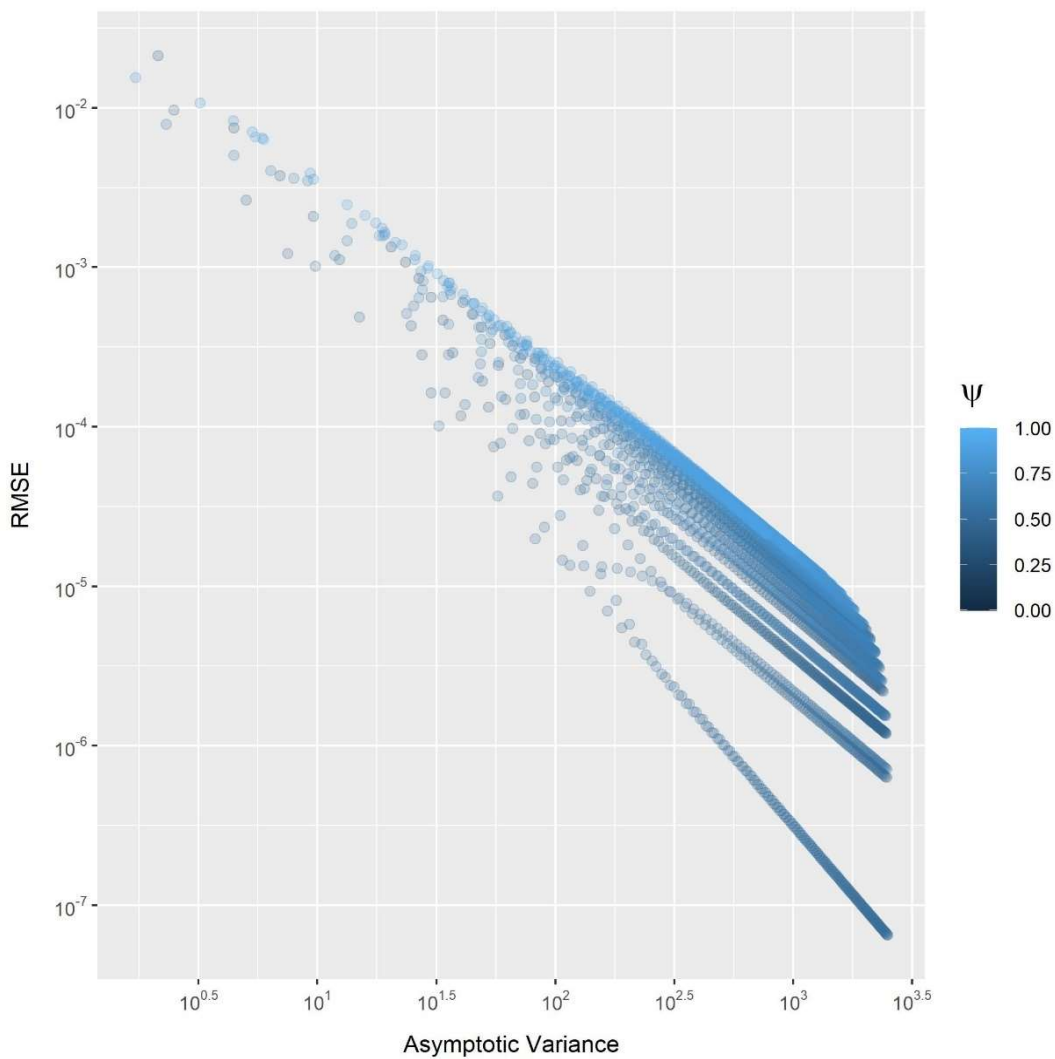

**FIGURE 3:** RMSE versus asymptotic variance of approximation

---

[6] There are more points in Figure 3 than there are squares in Figure 2. This is because we have limited the data for Figure 2 to show only cases up to $n = 1000$ for illustrate purposes.

[7] Running a linear regression on these simulation values with LRMSE as the response variable and the asymptotic log-variance and $\psi$ as the explanatory variables gives $R^2 = 0.9915$. There is little additional explanatory power added by adding $n$, $k$ or $\phi$.

We can see from Figures 2-3 that the present approximation to the spillage distribution works well when the parameters give a high variance for the distribution. In Figure 3 we have used the asymptotic variance $\sigma_{n,k,\phi}^{*2}$ and the proportion $\psi$, which can be computed non-recursively, but the general negative relationship holds when using the exact variance. The reader should bear in mind that our present accuracy results are only for the approximation in Bender (1973) and there are plenty of other approximations available. Consequently, our present analysis is not the last word in approximation accuracy for the spillage distribution, but it gives a sense of the accuracy over a range of parameter combinations for one of the possible approximation. What matters is that it is feasible to approximate the probability mass function of the spillage distribution using a non-recursive method that is computationally feasible for large $n$ and $k$, and this approximation performs well when the variance of the distribution becomes large.

The approximation used here (and taken from Bender 1973) is implemented in the probability functions **`dspillage`**, **`pspillage`** and **`qspillage`** in the **`occupancy`** package. The user can set **`approx = TRUE`** in the input to the functions to use the approximation instead of recursive computation. This gives a simple user-friendly facility to compute the probability functions for the spillage distribution either using the "exact" recursive method or the faster approximation method. As a general rule, we recommend that users compute the distribution using the recursive method when feasible but user the approximation when $n$ is large enough to give a high variance for the distribution.

**4. Concluding remarks**

In the present paper we have undertaken an examination of the spillage distribution arising in occupancy analysis in the extended balls-in-bins model. This probability distribution describes the behaviour of the spillage number conditional on the occupancy number. We have seen that the distribution and its moments involve the noncentral Stirling numbers of the second kind. We have shown how to compute the mass function for the distribution using recursive methods and how to compute an approximation that does not require recursion. We have also examined the accuracy of this approximation and we find that it works well in cases where the variance of the distribution is large.

Our present analysis extends an initial analysis of occupancy distributions in O'Neill (2021), giving more detail on the properties of the spillage distribution. We hope that this examination will pique the reader's interest in the spillage distribution and other distributions arising in the occupancy problems in this broader analysis.

# Appendix: Proof of Theorems

**PROOF OF THEOREM 1:** We begin by proving the form of the probability generating function for all $z > 0$. We can use our expansion for the noncentral Stirling numbers of the second kind to obtain:

$$G(z) = \mathbb{E}(z^{R_n}) = \sum_{r=0}^{n-k} z^r \cdot \text{Spillage}(r|n, k, \phi)$$

$$= \sum_{r=0}^{n-k} z^r \cdot \binom{n}{k+r} \cdot \phi^{n-k-r} \cdot \frac{S(k+r, k)}{S(n, k, \phi)}$$

$$= \frac{z^{n-k}}{S(n, k, \phi)} \sum_{r=0}^{n-k} \binom{n}{k+r} \cdot \left(\frac{\phi}{z}\right)^{n-k-r} \cdot S(k+r, k)$$

$$= z^{n-k} \cdot \frac{S(n, k, \phi/z)}{S(n, k, \phi)}.$$

The other generating functions follow using simple substitutions and transformations of the probability generating function. ∎

**PROOF OF LEMMA 1:** From O'Neill (2021) (Lemma 1) we know that the non-central Stirling numbers of the second kind obey the partial differential equation:

$$\frac{\partial S}{\partial \phi}(n, k, \phi) = n \cdot S(n-1, k, \phi).$$

Using the chain rule and applying this result gives:

$$\frac{dS}{ds}(n, k, \phi e^{-s}) = -\phi e^{-s} \cdot \frac{\partial S}{\partial \phi}(n, k, \phi e^{-s}) = -n\phi e^{-s} \cdot S(n-1, k, \phi).$$

Consequently, differentiating $Q_\ell$ and applying the quotient rule and then the chain rule gives:

$$\frac{dQ_\ell}{ds}(s) = \frac{d}{ds} \frac{S(n-\ell, k, \phi e^{-s})}{S(n, k, \phi e^{-s})}$$

$$= \frac{1}{S(n, k, \phi e^{-s})^2} \left[ \begin{array}{l} \frac{dS}{ds}(n-\ell, k, \phi e^{-s}) \cdot S(n, k, \phi e^{-s}) \\ -S(n-\ell, k, \phi e^{-s}) \cdot \frac{dS}{ds}(n, k, \phi e^{-s}) \end{array} \right]$$

$$= \frac{\phi e^{-s}}{S(n, k, \phi e^{-s})^2} \left[ \begin{array}{l} -(n-\ell) \cdot S(n-\ell-1, k, \phi e^{-s}) \cdot S(n, k, \phi e^{-s}) \\ +n \cdot S(n-\ell, k, \phi e^{-s}) \cdot S(n-1, k, \phi e^{-s}) \end{array} \right]$$

$$= \phi e^{-s} \left[ \begin{array}{c} n \cdot \dfrac{S(n-\ell, k, \phi e^{-s})}{S(n, k, \phi e^{-s})} \cdot \dfrac{S(n-1, k, \phi e^{-s})}{S(n, k, \phi e^{-s})} \\ -(n-\ell) \cdot \dfrac{S(n-\ell-1, k, \phi e^{-s})}{S(n, k, \phi e^{-s})} \end{array} \right]$$

$$= \phi e^{-s} [n Q_1(s) Q_\ell(s) - (n-\ell) Q_{\ell+1}(s)],$$

which was to be shown. ∎

**PROOF OF LEMMA 2:** This proof uses a rule in O'Neill (2021) (Lemma 3) for moving the non-centrality parameter in the noncentral Stirling numbers of the second kind:

$$S(n, k, \phi') = \sum_{\ell=0}^{n-k} \binom{n}{\ell} (\phi' - \phi)^\ell \cdot S(n-\ell, k, \phi).$$

Substituting $\phi' = \phi e^{-t}$ we get the alternative form:

$$S(n, k, \phi e^{-t}) = \sum_{\ell=0}^{n-k} \binom{n}{\ell} (e^{-t} - 1)^\ell \phi^\ell \cdot S(n-\ell, k, \phi).$$

We can therefore write the moment generating function of the distribution as:

$$m(t) = e^{t(n-k)} \cdot \frac{S(n, k, \phi e^{-t})}{S(n, k, \phi)}$$

$$= e^{t(n-k)} \sum_{\ell=0}^{n-k} \binom{n}{\ell} (e^{-t} - 1)^\ell \cdot \phi^\ell \cdot \frac{S(n-\ell, k, \phi)}{S(n, k, \phi)}$$

$$= e^{t(n-k)} \sum_{\ell=0}^{n-k} \binom{n}{\ell} (e^{-t} - 1)^\ell \cdot H_\ell,$$

which was to be shown. ∎

**PROOF OF THEOREM 2:** Here we will derive the moments in the theorem from the probability generating function shown in Theorem 1. Using the differentiation result in Lemma 1, the first four derivatives of the cumulant generating function are:

$$\frac{dK}{ds}(s) = (n - k) - n\phi e^{-s} \cdot Q_1(s),$$

$$\frac{d^2 K}{ds^2}(s) = n\phi e^{-s} \cdot Q_1(s) - n^2 \phi^2 e^{-2s} \cdot Q_1(s)^2$$
$$+ n(n-1)\phi^2 e^{-2s} \cdot Q_2(s),$$

$$\frac{d^3K}{ds^3}(s) = -n\phi e^{-s} \cdot Q_1(s) + 3n^2\phi^2 e^{-2s} \cdot Q_1(s)^2$$
$$- 3n(n-1)\phi^2 e^{-2s} \cdot Q_2(s) - 2n^3\phi^3 e^{-3s} \cdot Q_1(s)^3$$
$$+ 3n^2(n-1)\phi^3 e^{-3s} \cdot Q_1(s)Q_2(s)$$
$$- n(n-1)(n-2)\phi^3 e^{-3s} \cdot Q_3(s),$$

$$\frac{d^4K}{ds^4}(s) = n\phi e^{-s} \cdot Q_1(s) - 7n^2\phi^2 e^{-2s} \cdot Q_1(s)^2$$
$$+ 7n(n-1)\phi^2 e^{-2s} \cdot Q_2(s) + 12n^3\phi^3 e^{-3s} \cdot Q_1(s)^3$$
$$- 18n^2(n-1)\phi^3 e^{-3s} \cdot Q_1(s)Q_2(s) + 6n(n-1)(n-2)\phi^3 e^{-3s} \cdot Q_3(s)$$
$$- 6n^4\phi^4 e^{-4s} \cdot Q_1(s)^4 + 12n^3(n-1)\phi^4 e^{-4s} \cdot Q_1(s)^2 Q_2(s)$$
$$- 3n^2(n-1)^2\phi^4 e^{-4s} \cdot Q_2(s)^2 - 4n^2(n-1)(n-2)\phi^4 e^{-4s} \cdot Q_1(s)Q_3(s)$$
$$+ n(n-1)(n-2)(n-3)\phi^4 e^{-4s} \cdot Q_4(s).$$

Substituting $s = 0$ and using the substitution $H_\ell = \phi^\ell \cdot Q_\ell(0)$ gives the cumulants:

$$\kappa_1 \equiv \frac{dK}{ds}(0) = (n-k) - nH_1,$$

$$\kappa_2 \equiv \frac{d^2K}{ds^2}(0) = nH_1 - n^2 H_1^2 + n(n-1)H_2,$$

$$\kappa_3 \equiv \frac{d^3K}{ds^3}(0) = -nH_1 + 3n^2 H_1^2 - 3n(n-1)H_2 - 2n^3 H_1^3$$
$$+ 3n^2(n-1)H_1 H_2 - n(n-1)(n-2)H_3,$$

$$\kappa_4 \equiv \frac{d^4K}{ds^4}(0) = nH_1 - 7n^2 H_1^2 + 7n(n-1)H_2 + 12n^3 H_1^3$$
$$- 18n^2(n-1)H_1 H_2 + 6n(n-1)(n-2)H_3$$
$$- 6n^4 H_1^4 + 12n^3(n-1)H_1^2 H_2 - 3n^2(n-1)^2 H_2^2$$
$$- 4n^2(n-1)(n-2)H_1 H_3 + n(n-1)(n-2)(n-3)H_4.$$

Converting the cumulants to the central moments then gives the formulae in the theorem. ∎

**PROOF OF THEOREM 3:** We will use the following recursive rule for the noncentral Stirling numbers of the second kind (see e.g., O'Neill 2021, Lemma 1):

$$S(n, k, \phi) = (k + \phi) \cdot S(n-1, k, \phi) + S(n-1, k-1, \phi).$$

Using repeated application of this rule and re-arranging gives:

$$S(n - \ell, k, \phi) = \frac{S(n, k, \phi)}{(k + \phi)^\ell} - \sum_{r=1}^{\ell} (k + \phi)^{r-\ell-1} \cdot S(n - r, k - 1, \phi).$$

This allows us to write $H_\ell$ as:

$$H_\ell = \phi^\ell \cdot \frac{S(n-\ell, k, \phi)}{S(n, k, \phi)}$$

$$= \frac{\phi^\ell}{S(n, k, \phi)} \left[ \frac{S(n, k, \phi)}{(k+\phi)^\ell} - \sum_{r=1}^{\ell} (k+\phi)^{r-\ell-1} \cdot S(n-r, k-1, \phi) \right]$$

$$= \left(\frac{\phi}{k+\phi}\right)^\ell - \phi^\ell \sum_{r=1}^{\ell} (k+\phi)^{r-\ell-1} \cdot \frac{S(n-r, k-1, \phi)}{S(n, k, \phi)}.$$

Consequently, we have:

$$\lim_{n \to \infty} H_\ell = \left(\frac{\phi}{k+\phi}\right)^\ell - \phi^\ell \sum_{r=1}^{\ell} (k+\phi)^{r-\ell-1} \cdot \lim_{n \to \infty} \frac{S(n-r, k-1, \phi)}{S(n, k, \phi)}$$

$$= \left(\frac{\phi}{k+\phi}\right)^\ell - \phi^\ell \sum_{r=1}^{\ell} (k+\phi)^{r-\ell-1} \cdot 0$$

$$= \left(\frac{\phi}{k+\phi}\right)^\ell.$$

This establishes the limit for $H_\ell$ in the theorem. Substituting the limiting form into the moment equations in Theorem 2 and using the asymptotic substitution $n \sim (n-k)$ then gives the asymptotic equivalence:

$$\mu_{n,k,\phi} \sim (n-k) - n\varphi$$
$$\sim (n-k) - (n-k)\varphi$$
$$= (n-k)(1-\varphi),$$

$$\sigma^2_{n,k,\phi} \sim n\varphi - n^2\varphi^2 + n(n-1)\varphi^2$$
$$= n\varphi(1-\varphi)$$
$$\sim (n-k)\varphi(1-\varphi),$$

$$\gamma_{n,k,\phi} \sim \frac{\begin{bmatrix} -n\varphi + 3n^2\varphi^2 - 3n(n-1)\varphi^2 - 2n^3\varphi^3 \\ +3n^2(n-1)\varphi^3 - n(n-1)(n-2)\varphi^3 \end{bmatrix}}{[n\varphi(1-\varphi)]^{3/2}}$$

$$= \frac{\begin{bmatrix} -n\varphi + 3n\varphi^2 - 2n^3\varphi^3 \\ +2n(n-1)(n+1)\varphi^3 \end{bmatrix}}{[n\varphi(1-\varphi)]^{3/2}}$$

$$= \frac{[-n\varphi + 3n\varphi^2 - 2n\varphi^3]}{[n\varphi(1-\varphi)]^{3/2}}$$

$$= \frac{n\varphi(2\varphi - 1)(1-\varphi)}{[n\varphi(1-\varphi)]^{3/2}}$$

$$= \frac{2(\varphi - \frac{1}{2})}{\sqrt{n\varphi(1-\varphi)}}$$

$$\sim \frac{2(\varphi - \frac{1}{2})}{\sqrt{(n-k)\varphi(1-\varphi)}},$$

$$\kappa_{n,k,\phi} \sim 3 + \frac{\begin{bmatrix} n\varphi - 7n^2\varphi^2 + 7n(n-1)\varphi^2 + 12n^3\varphi^3 \\ -18n^2(n-1)\varphi^3 + 6n(n-1)(n-2)\varphi^3 \\ -6n^4\varphi^4 + 12n^3(n-1)\varphi^4 - 3n^2(n-1)^2\varphi^4 \\ -4n^2(n-1)(n-2)\varphi^4 \\ +n(n-1)(n-2)(n-3)\varphi^4 \end{bmatrix}}{[n\varphi(1-\varphi)]^2}$$

$$= 3 + \frac{[n\varphi - 7n\varphi^2 + 12n\varphi^3 - 6n\varphi^4]}{[n\varphi(1-\varphi)]^2}$$

$$= 3 + \frac{n\varphi(1-\varphi)[1 - 6\varphi + 6\varphi^2]}{[n\varphi(1-\varphi)]^2}$$

$$= 3 + \frac{1 - 6\varphi + 6\varphi^2}{n\varphi(1-\varphi)}$$

$$= 3 + \frac{1 - 6\varphi(1-\varphi)}{n\varphi(1-\varphi)}$$

$$= 3 + \frac{1 - 6\varphi(1-\varphi)}{(n-k)\varphi(1-\varphi)}.$$

This establishes the asymptotic forms in the theorem. ∎

**PROOF OF THEOREM 4:** The probability mass function for the occupancy distribution is:

$$\text{Occ}(k|n, m, \theta) = \frac{\theta^n}{m^n} \cdot (m)_k \cdot S\left(n, k, m \cdot \frac{1-\theta}{\theta}\right),$$

and the relevant spillage probabilities are:

$$\text{Spillage}\left(x - k \middle| n, k, m \cdot \frac{1-\theta}{\theta}\right) = \binom{n}{x} \cdot \left(m \cdot \frac{1-\theta}{\theta}\right)^{n-x} \cdot S(x, k) \Big/ S\left(n, k, m \cdot \frac{1-\theta}{\theta}\right).$$

We therefore have:

$$\text{RHS} = \sum_{k=0}^{x} \text{Spillage}\left(x - k \middle| n, k, m \cdot \frac{1-\theta}{\theta}\right) \cdot \text{Occ}(k|n, m, \theta)$$

$$= \sum_{k=0}^{x} \binom{n}{x} \cdot \left(m \cdot \frac{1-\theta}{\theta}\right)^{n-x} \cdot S(x, k) \cdot \frac{\theta^n}{m^n} \cdot (m)_k$$

$$= \binom{n}{x} \cdot \theta^x (1-\theta)^{n-x} \sum_{k=0}^{x} S(x, k) \cdot \frac{(m)_k}{m^x}$$

$$= \binom{n}{x} \cdot \theta^x (1-\theta)^{n-x} = \text{Bin}(x|n,\theta),$$

which was to be shown. ∎

**PROOF OF THEOREM 5:** The base equation for the recursion follows directly from the fact that $S(k,k) = 1$ for all $k \geq 0$. To prove the recursive equation we define:

$$K(r,k) = \binom{n}{k+r} \cdot \phi^{n-k-r} \cdot S(k+r, k),$$

The Stirling numbers obey the recursive equation $S(n+1, k) = k \cdot S(n,k) + S(n, k-1)$ so we can apply this to obtain:

$$\begin{aligned}
K(r+1, k) &= \binom{n}{k+r+1} \cdot \phi^{n-k-r-1} \cdot S(k+r+1, k) \\
&= \binom{n}{k+r+1} \cdot \phi^{n-k-r-1} \cdot [k \cdot S(k+r, k) + S(k+r, k-1)] \\
&= \frac{1}{\phi} \cdot \frac{n-k-r}{k+r+1} \cdot \binom{n}{k+r} \cdot \phi^{n-k-r} \cdot [k \cdot S(k+r, k) + S(k+r, k-1)] \\
&= \frac{n-k-r}{k+r+1} \cdot \frac{k \cdot K(r, k) + K(r+1, k-1)}{\phi}.
\end{aligned}$$

Taking logarithms gives the corresponding recursive equation in the theorem. The final result follows from the expansion from the noncentral Stirling numbers. ∎